\documentclass[10pt]{amsart}
\usepackage{graphicx}                                    
\usepackage{amssymb}
\usepackage{enumitem}
\usepackage{color}
\usepackage[all]{xy}
\usepackage{comment}

\newtheorem{theorem}{Theorem}[section] 
\newtheorem{lemma}[theorem]{Lemma}
\newtheorem{conjecture*}[theorem]{Conjecture}
\newtheorem{proposition}[theorem]{Proposition}
\newtheorem{conjecture}[theorem]{Conjecture}

\newtheorem{corollary}[theorem]{Corollary}
\newtheorem{remark}[theorem]{Remark}
\newtheoremstyle{notauto}{}{}{\itshape}{}{\bfseries}{.}{0.5em}{\thmnote{#3}}
\theoremstyle{notauto}

\theoremstyle{definition}
\newtheorem{definition}[theorem]{Definition}
\newtheorem{example}[theorem]{Example}
\theoremstyle{remark}

\newcommand{\bbF}{{\mathbb F}}
\newcommand{\bbZ}{{\mathbb Z}}
\newcommand{\bbR}{{\mathbb R}}

\newcommand{\res}{\mathrm{res}}
\newcommand{\rk}{\mathrm{rk}}
\newcommand{\im}{\mathrm{Im}}

\begin{document}
\baselineskip13pt
	
\title[Free actions on 
products of real projective spaces]{Free actions on products of real projective spaces}

\author{Erg{\" u}n Yal{\c c}{\i}n}
\address{Department of Mathematics, Bilkent University, 06800 
Bilkent, Ankara, Turkey}

\email{yalcine@fen.bilkent.edu.tr}

\begin{abstract}  
We prove that if $G=(\mathbb{Z}/2)^r$ acts freely and cellularly on a finite-dimensional CW-complex $X$ homotopy equivalent to $\mathbb{R}P ^{n_1} \times \cdots \times \mathbb{R} P ^{n_k}$ with trivial action on the mod-$2$ cohomology, then $r \leq \mu (n_1)+ \cdots + \mu(n_k )$
where for each integer $n\geq 1$,  $\mu (n)=0$ if $n$ is even, 
$\mu(n)=1$ if $n\equiv 1 \pmod 4$, and $\mu(n)=2$ if $n\equiv 3 \pmod 4$.
This proves a homotopy-theoretic version of a conjecture of Cusick.
\end{abstract}

\date{August 20, 2026}



\maketitle

\section{Introduction}

Classifying all finite groups that can act freely on a product of spheres $X=S^{n_1} \times \cdots \times S^{n_k}$ is one of the most important open problems in transformation group theory.  There is a conjecture that states  that if $G=(\mathbb{Z}/p)^r$ acts freely on a product of $k$ spheres, then $r\leq k$.  Special cases of this conjecture were proved by Carlsson \cite{Carlsson-1980}, \cite{Carlsson-1982}, Adem and Browder \cite{AdBr-1988}, Yal\c c\i n \cite{Yalcin-2000}, Hanke \cite{Hanke-2009}, Okutan and Yal\c c\i n \cite{OkYa-2013}, and others.  

The techniques developed by Carlsson in \cite{Carlsson-1980} and \cite{Carlsson-1982} 
for studying free actions on products of spheres were used later to study 
free actions on other spaces such as 
products of real projective spaces \cite{Cusick-1983}, products of Milnor manifolds \cite{Singh-2014}, and 
products of Dold spaces 
\cite{Dey-2023}. In the 1980s, a sequence of papers was published by Cusick on free actions on products of even-dimensional spheres and on products of real projective spaces (see  \cite{Cusick-1983},  \cite{Cusick-1987}).
In one of these papers Cusick stated the following  conjecture:

\begin{conjecture}[Cusick \cite{Cusick-1983}] \label{conj:Cusick} 
Suppose that $X$ is a finite CW-complex with the mod-$2$ cohomology of  $ \prod_{i=1} ^k \mathbb{R}P ^{n_i}$.
If $X$ admits a free action of  $G=(\bbZ/2)^r$ in such a way that the induced action on mod-$2$ cohomology is trivial,
then 
\[
r \leq \mu (n_1)+ \cdots + \mu(n_k ),
\]
where for each integer $n\geq 1$, 
\begin{equation*}
\mu(n) = \begin{cases} 0 & \text{if } \ n \text{ is even}, \\
1 & \text{if } \ n\equiv 1 \pmod 4, \\
2 & \text{if } \  n\equiv 3 \pmod 4.
\end{cases}
\end{equation*} 
\end{conjecture}

Cusick \cite{Cusick-1983} proved Conjecture \ref{conj:Cusick} in the case where $n_i \not \equiv 3 \pmod 4$ for all $i=1,\dots, k$.
Later Yal{\c c}{\i}n \cite{Yalcin-2000} proved the conjecture for the case where $X$ is homotopy equivalent to  
$ \prod_{i=1} ^k \mathbb{R}P ^{n_i}$ where all $n_i$ are odd, under the additional assumption that the induced action on integral cohomology is trivial (see Remark \ref{rem:Mistake}).  The equidimensional case without the assumption that the action on the mod-$2$ cohomology is trivial  was settled
by Adem and Yal{\c c}{\i}n  \cite{AdYa-1999}. In this paper we prove the following:

\begin{theorem}\label{thm:Main}  Suppose that $G=(\mathbb{Z}/2)^r$ and $X$ is a finite-dimensional free $G$-CW-complex 
homotopy equivalent to $\prod _{i=1} ^k  \mathbb{R}P ^{n_i}$ 
such that the $G$-action on the mod-$2$ cohomology of $X$ is trivial. Then 
\[
r \leq \mu (n_1)+ \cdots + \mu(n_k ),
\]
where for each integer $n\geq 1$, 
\begin{equation*}
\mu(n) = \begin{cases} 0 & \text{if } \ n \text{ is even}, \\
1 & \text{if } \ n\equiv 1 \pmod 4, \\
2 & \text{if } \  n\equiv 3 \pmod 4.
\end{cases}
\end{equation*} 
\end{theorem}

Note that the assumption that $X$ is homotopy equivalent to a product of real projective spaces is stronger than the assumption in Cusick's conjecture. Since we use integral cohomology and integral Bockstein operator, an assumption stronger than the mod-$2$ cohomology isomorphism is necessary for our proof. On the other hand, the assumption that $X$ is finite-dimensional is weaker than that $X$ is a finite CW-complex.

There are examples of free actions that show that the upper bound given in Conjecture \ref{conj:Cusick} is sharp. To construct such actions, first observe that the quaternion group $Q_8$ of order $8$ acts freely on $S^3$ via its identification as the group of unit quaternions. By taking joins of this action, we obtain a $Q_8$-action on $S^{4m+3}$ for all $m\geq 0$ which induces a free $\bbZ/2\times \bbZ/2$-action on $\bbR P^{4m+3}$ for all $m\geq 0$. There is also a free $\bbZ /2$-action on $\bbR P^{2m+1}$ for any $m\geq 0$, induced by the free $\bbZ/4$-action on $S^{2m+1}$ defined by $(z_0, \dots, z_m)\to (iz_0, \dots, iz_m)$. 
Taking products of these actions we obtain a free $G=(\bbZ/2)^{k+2l}$ action on $X=\prod _{i=1} ^{k} \bbR P^{2m_i+1} \times \prod _{i=1}^{l} \bbR P^{4m_i+3}$ with trivial action on the integral cohomology (see Example \ref{ex:Trivial} for details). 

There are also examples of free actions on products of real projective spaces that induce a nontrivial action on integral cohomology and trivial action on mod-$2$ cohomology. To construct such actions, consider the $D_8$-action on $S^{2m+1}$ obtained by taking joins of the $D_8$-action on the circle $S^1$ defined by the symmetries of a square. This action induces a $\bbZ/2 \times \bbZ/2$-action on $\bbR P^{2m+1}$ with nontrivial action on the integral cohomology when $m$ is even. Combining this action with another action obtained via inflation from a free $\bbZ/2$-action on $\bbR P ^{2m+1}$, and taking products of these actions, we can obtain a free $(\bbZ/2)^{k+l}$-action
on $\prod _{i=1}^k  \bbR P^{2m_i+1} \times \prod _{i=1}^l \bbR P ^{4m_i+1}$ for any $k\geq l\geq 1$ with nontrivial action on integral cohomology whose mod-$2$ reduction is trivial (see Examples \ref{ex:Nontrivial} and \ref{ex:NontrivialFree}). 

We prove Theorem \ref{thm:Main} using the methods from \cite{Cusick-1983} and \cite{Yalcin-2000}. 
The new input comes from a closer look at the interaction between the Serre spectral sequences with $\bbF_2$- and $\bbZ$-coefficients. Recall that if $X$ is a $G$-CW-complex, there is an associated Borel fibration $X \to X_G \to BG$ where $X_G=(EG\times X)/G$.
For a commutative ring $R$, the Serre spectral sequence for the Borel fibration with coefficients in $R$ is of the form
\[
E_2 ^{p,q} =H^p (BG; H^q (X; R) ) \Rightarrow H^{p+q} (X_G; R).
\]

The Serre spectral sequence has a multiplicative structure and if $R=\bbF_2$ and the induced $G$-action on $H^* (X; \bbF_2)$  is trivial, then there is an isomorphism of bigraded $\bbF_2$-algebras 
\[
E_2 ^{*, *} \cong H^* (BG; \bbF_2) \otimes H^* (X; \bbF_2).
\]
If $G\cong (\bbZ/2)^r$, then we have   
\[
H^* (G; \bbF_2) \cong \bbF_2 [x_1, \dots, x_r]
\]
where $x_i=\pi_i^*(x)$, with $x$ the generator of $H^1(\bbZ/2;\mathbb F_2)$. Similarly, for $X\simeq \prod _{i=1} ^k \mathbb{R} P^{n_i}$, we have  
\[
H^* (X; \bbF_2) \cong \bbF_2 [t_1, \dots, t_k]/ ( t_1 ^{n_1+1}, \dots, t_k ^{n_k+1})
\]
where $t_i=\pi_i^*(t)$, with $t$ the generator of $H^1(\mathbb RP^{n_i};\mathbb F_2)$.

Let  $d_2: E_2 ^{0, 1} \to E_2 ^{2, 0}$ be the differential on the $E_2$-page, and $\alpha_1, \dots , \alpha_k \in H^2 (G; \bbF_2) $ be the 2-dimensional cohomology classes such that 
$d_2 (1 \otimes t_i)=\alpha_i \otimes 1$ for each $i$. We call these cohomology classes the k-invariants of the $G$-action  on $X$, a terminology that goes back to the obstruction theory.

It is proved by Cusick \cite{Cusick-1983} that for every $i$ where $n_i$ is even, $\alpha_i=0$, and  for every $i$ where $n_i\equiv 1 \pmod 4$, $Sq^1 (\alpha _i)$ lies in the ideal generated by 
$\alpha_1, \dots, \alpha_k$ in $H^* (G; \bbF_2)$. Using these Cusick proved the conjecture under the assumption that for all $i$, $n_i \not \equiv 3 \pmod 4$ (see \cite[Main Theorem]{Cusick-1983}).

Using a result by Carlsson \cite{Carlsson-1980}, it can be shown that if the $G$-action on $X$ is free, then the k-invariants 
$\alpha_1, \dots, \alpha_k$ have no nontrivial common zeros when they are considered as quadratic polynomials in $\bbF_2$-coefficients (see Proposition \ref{pro:NoCommonZeros}). This gives the inequality also for the case where $n_i \not \equiv 1 \pmod 4$ for all $i$ (see \cite[Main Theorem]{Cusick-1983}). 

To prove the inequality in Theorem \ref{thm:Main} under the assumption that $X$ is homotopy equivalent to $\prod _{i=1} ^k \bbR P^{n_i}$ with trivial action on mod-2 cohomology, we first consider the case where the induced $G$-action on the integral cohomology of $X$ is trivial. Under the assumption that the $\bbF_2$-vector space generated by the $k$-invariants $\alpha_1,\dots,\alpha_k$ has an irreducible basis (see Definition \ref{def:Irreducible}), we prove that for all the indices $i$ of the basis elements, $n_i\equiv 3 \pmod 4$ (see Lemma \ref{lem:Square}). Together with the induction argument used in the proof of Theorem \ref{thm:Main}, this proves the desired inequality in the case where the integral cohomology action is trivial.

In Section \ref{sect:Nontrivial}, we consider the free actions with nontrivial induced $G$-action on the integral cohomology of $X$. We first give an example of such action in Example \ref{ex:NontrivialFree}, and then analyze the integral Serre spectral sequence in this case. We prove that under the assumption that the $\bbF_2$-vector space generated by $k$-invariants has an irreducible basis, all indices $i$ of the basis elements satisfy $n_i\equiv 3 \pmod 4$ (see Proposition \ref{pro:Main}). We show in Section \ref{sect:MainThm} that this fact proves the inequality in Theorem \ref{thm:Main} by induction.

\subsection{Acknowledgements} The author is supported by T\" UB\. ITAK grant 2219-International Postdoctoral Research Fellowship Program (2023, 2nd term). We thank T\" UB\. ITAK for their support of this research. We also thank Matthew Gelvin and Caroline Yal\c c\i n for reading an earlier version of the paper and for their suggestions. 
We thank Marc~Stephan for providing the reference for Proposition 9.20 in \cite{Osborne-2000}. We also thank the referee for careful reading of the paper and suggestions that improved the paper.


\section{The Serre spectral sequence for the Borel construction}\label{sect:Serre}

Let $G$ be a finite group and $X$ a $G$-CW-complex. The Borel construction
$X_G=EG\times_G X$ fits into the Borel fibration
\[
X\to X_G\to BG
\]
(see \cite[Proposition II.2.4]{Bredon-1972}). For a commutative ring \(R\), the associated Serre spectral sequence is of the form
\[
E_2^{p,q}=H^p(BG;\mathcal H^q(X;R))\Rightarrow H^{p+q}(X_G;R),
\]
where the local system is determined by the induced $G$-action on $H^q(X;R)$ (see \cite[Theorem 5.2]{McCleary-2001}).  

When the $G$-action on the cohomology of $X$ is trivial, the local coefficient system is trivial
and $E_2 ^{p, q}$ is isomorphic to $H^p (BG; H^q (X;R))$.  In this case we say the Borel fibration is simple. If $R$ is a field and the Borel fibration is simple, then there is an isomorphism of $R$-algebras 
\[
E_2 ^{p,q} =H^p (BG; H^q (X; R))\cong 
H^p (BG; R) \otimes _R H^q (X; R).
\]

The Serre spectral sequence has a product structure defined by the composition 
\[
H^p (BG; \mathcal{H} ^q( X; R)) \otimes H^s(BG; \mathcal{H}^t (X; R)) \to H^{p+s} (BG; \mathcal{H} ^q (X; R) \otimes \mathcal{H}^t (X; R)) 
\]
\[
\to H^{p+s} (BG; \mathcal{H} ^{q+t} (X; R)),
\]
where both maps are defined by the cup product.
In the case where the Borel fibration is simple and $R$ is a field, the product structure on $E_2 ^{*,*}\cong H^p (BG; R) \otimes _R H^q (X; R)$ is defined by $(x \otimes y) (u\otimes v) = xu \otimes yv$, where $xu$ and $yv$ are the products defined by the cup products in $H^* (BG; R)$ and $H^* (X; R)$ (see \cite[Theorem 5.6]{McCleary-2001}).

The Serre spectral sequence of the Borel fibration is useful for finding obstructions to the existence of free actions on finite-dimensional spaces because of the following well-known observation:

\begin{lemma}\label{lem:E_infty} Let $G$ be a finite group and $X$ a finite-dimensional free $G$-CW-complex. If $E_r ^{*, *}$ denotes the $r$-th page of the Serre
spectral sequence for the Borel fibration $X_G \to BG$ with $\bbF_2$-coefficients, then
$E_{\infty} ^{p,q}=0$  for all $p,q$ with $p+q> \dim X$.
\end{lemma}

\begin{proof} If the $G$-action on $X$ is free, then $X_G = EG\times _G X$ is homotopy equivalent to the orbit space $X/G$ 
(see \cite[Proposition 1]{Carlsson-1980}). Since $X$ is a finite-dimensional $G$-CW-complex, 
$X/G$ is also finite-dimensional with dimension equal to the dimension of $X$. This gives $H^n (X_G; \bbF_2)=0$ for all $n>\dim X$, hence $E^{p,q}_{\infty} =0$ for all $p,q$ with $p+q>\dim X$.
 \end{proof}

Let $G=(\bbZ/2)^r$ and $X$ be a finite-dimensional $G$-CW-complex
whose mod-$2$ cohomology  is isomorphic to the cohomology of $\prod _{i=1} ^k \mathbb{R}P ^{n_i}$.
Suppose that $G$ acts freely on $X$ in such a way that the induced $G$-action  
on $H^* (X; \bbF_2)$ is trivial.  Then there is an isomorphism of bigraded $\bbF_2$-algebras
\[
E_2 ^{*, *} \cong H^* (BG; \bbF_2) \otimes H^* (X; \bbF_2),
\]
where 
\[
H^* ( BG; \bbF_2)  \cong \bbF_2 [x_1, \dots, x_r]
\]
and 
\[
H^* (X; \bbF_2)  \cong \bbF_2 [t_1, \dots, t_k] / ( t_1 ^{n_1+1}, \dots, t_k ^{n_k+1} )
\]
with  $|x_i|=1$ for all $i=1,\dots, r$, and $|t_i|=1$ for all $i=1,\dots, k$.

\begin{definition} Let $\alpha_1, \dots, \alpha_k \in H^2 (G; \bbF_2)$ be the cohomology classes such that
 $d_2 (1 \otimes t_i)=\alpha_i \otimes 1$  for $i=1, \dots, k$.  
 These cohomology classes are called the \emph{k-invariants} of the $G$-action on $X$.
\end{definition}

When it is convenient, we identify the vertical line $E_2 ^{0, *} \cong  \bbF_2 \otimes H^* (X; \bbF_2)$ with $H^* (X;\bbF_2)$ and the horizontal line $E_2 ^{*, 0} \cong H^* (BG; \bbF_2) \otimes \bbF_2$ with $H^* (BG; \bbF_2)$ and write the elements in $E_2 ^{p, q}$ as products $uv$ where $u\in H^p (BG; \bbF_2)$ and $v\in H^q (X; \bbF_2)$. With these identifications, we sometimes write $d_2(t_i)=\alpha_i$.

In \cite{Cusick-1983}, Cusick proves the following:

\begin{lemma}[{\cite[Proposition B]{Cusick-1983}}]\label{lem:Properties}  (a) For every $i$ where $n_i$ is even,  $\alpha_i=0$. \\
(b) For every $i$ where $n_i \equiv 1 \pmod 4$, the Steenrod power $Sq ^1 (\alpha_i) $ lies in the ideal generated by $\alpha_1, \dots, \alpha_k$ in $H^*(G; \bbF_2)$.
\end{lemma}

\begin{proof} Since the proofs are short, we repeat them here for the convenience of the reader. 

(a) Let $n_i=2m_i$ for some integer $m_i \geq 0$. Then $t_i ^{2m_i +1}=0$. This gives 
\[
0=d_2 (t_i ^{2m_i +1})= d_2 (t_i) t_i ^{2m_i}=\alpha_i t_i ^{2m_i}.
\]
On the $E_2$-page, multiplication by $t_i ^{2m_i}$ is injective. Hence $\alpha_i=0$.

(b) Let $n_i =4k_i+1$ for some $k_i \geq 0$, and let $u_i:=t_i^2$. Since $d_2 (u_i)=d_2(t_i^2)=0$, $u_i $ survives to the $E_3$-page.
The differential $d_2 : E_2 ^{0, 1} \to E_2 ^{2,0} $ is a transgression (see Definition \ref{def:Transgression}), so it commutes with 
Steenrod operations (see \cite[Theorem 6.8 and Corollary 6.9]{McCleary-2001}). This gives \[
d_3 (u_i)=  d_3 (Sq^1 (t_i))= Sq^1 (d_2 (t_i))=Sq^1 (\alpha_i).
\]

Since $u_i ^{2k_i+1}=0$, we have 
\[
0=d_3 (u_i ^{2k_i +1} )= d_3 (u_i) u_i ^{2k_i}= Sq^1 (\alpha _i ) u_i ^{2k_i}
\]
in $E_3 ^{3, 4k_i}$. Hence there is an element $\theta \in E_2^{1, 4k_i+1}$ such that $d_2 (\theta)=Sq^1 ( \alpha_i) u_i ^{2k_i}$. From the product structure of the $E_2$-page, we see that $d_2(\theta)$ must be equal to $d_2(f) t_i ^{4k_i}$
for some $f\in E_2 ^{1, 1}$. 
Since multiplication with $u_i ^{2k_i}$ in injective on the $E_2$-page,
we obtain $Sq^1 (\alpha_i)=d_2(f)$. If $f=\sum _j l_j t_j$ where $l_j\in H^1 (G, \bbF_2)$, then $Sq^1 (\alpha_i)=d_2 (f)=\sum _j l_j \alpha_j$. Hence $Sq^1 (\alpha_i )$ lies in the ideal generated by $\alpha_1, \dots, \alpha_k$.
\end{proof}

We also have the following observation for k-invariants. 

\begin{proposition}\label{pro:NoCommonZeros}
Let $G=(\bbZ/2)^r$ and $X$ be a finite-dimensional $G$-CW-complex
whose mod-$2$ cohomology ring is isomorphic to $ H^* (\prod _{i=1} ^k \mathbb{R}P ^{n_i} ; \bbF_2).$
Suppose that $G$ acts freely on $X$ and the induced action  
on $H^* (X; \bbF_2)$ is trivial. Then the k-invariants $\alpha_1, \dots, \alpha_k \in H^2 (G; \bbF_2)$ 
have no common zeros in $(\bbF_2)^r\backslash \{(0,\dots, 0)\}$ when they are considered as quadratic polynomials in the variables $x_1, \dots, x_r$.
\end{proposition}

\begin{proof} We repeat the argument given for the proof of the main theorem in \cite{Cusick-1983}. A similar argument for free actions on products of spheres was given by Carlsson in \cite{Carlsson-1980}.

Let $\alpha_1, \dots, \alpha_k \in H^2 (G; \bbF_2)$ be the k-invariants of the action. Consider the k-invariants as quadratic polynomials in the variables $x_1, \dots, x_r$. Let $g_1, \dots, g_r$ be the basis for $G=(\bbZ/ 2)^r$ dual to the basis $x_1, \dots, x_r$. Let $C\leq G$ be a nontrivial cyclic subgroup
of $G$ generated by $g= g_1 ^{\lambda_1}\dots g_r ^{\lambda _r }$, and $x\in H^1 (C; \bbF_2)$ be the generator dual to $g$.
Then $\res ^G _C \alpha _i $ is equal to the class $\alpha_i (\lambda _1, \dots, \lambda _r ) \cdot x^2$.  
Hence $\alpha_1, \dots, \alpha_k$ have no common zeros in $(\bbF_2)^r\backslash \{ (0, \dots, 0)\}$ as quadratic polynomials if and only if for every nontrivial cyclic subgroup 
$C \leq G$,  there is an $i \in \{1, \dots, k\}$ such that $\res ^G _C \alpha _i \neq 0$.

Assume, to the contrary, that there exists a nontrivial cyclic subgroup $C\leq G$ such that $\res^G _C \alpha _i=0$ for all $i$. 
Consider the restriction of the 
$G$-action on $X$ to the subgroup $C$. By the naturality of the Serre spectral sequence with respect to subgroup inclusions (see \cite[page 1150]{Carlsson-1980}),  in the Serre spectral sequence for $X_C \to BC$ we have $d_2 ( t_i )=\res ^G _C \alpha_i=0$ for all $i$.
By the product structure of the Serre spectral sequence, this implies that all the other differentials are
also zero. Hence the Serre spectral sequence for the $C$-action on $X$ collapses at the $E_2$-page.
This implies that $E_{\infty} ^{p, 0} = E_2 ^{p, 0}  \cong H^p(C; \bbF_2)\otimes H^0 (X; \bbF_2)$ is nonzero for all $p\geq 0$, which contradicts the conclusion of Lemma \ref{lem:E_infty}.
\end{proof}


\section{Comparison with the integral Serre spectral sequence}\label{sect:IntegralSS}

Let $F \xrightarrow{i}  E \xrightarrow{\pi} B$ be a fibration  where both $B$ and $F$ are path-connected. The Serre spectral sequence for the fibration $E \xrightarrow{\pi} B$ with coefficients in the abelian group $A$ is of the form
\[
E_2 ^{p,q} (A) = H^p (B; H^q (F; A)) \Rightarrow H^{p+q} (E; A).
\]
We denote the $r$-th page of this spectral sequence by $E_r ^{*, *} (A)$.

For every homomorphism $f: A \to A'$ between two abelian groups, there is a map of spectral sequences
\[
f_*: E_r ^{*,*} (A) \to E_r ^{*, *} (A').
\]
By this we mean that $f_* : E_r ^{*,*} (A) \to E_r ^{*,*} (A')$ is a homomorphism of bigraded abelian groups that commutes with the differentials $d_r ^{*, *} (A)$ and $d_r ^{*,*} (A')$ where 
$d_r ^{p,q} (A) : E_r ^{p,q} (A) \to E_r ^{p+r, q-r+1} (A)$.

The vertical edge of the Serre spectral sequence is $E_2 ^{0, q}\cong H^0 (B; H^q (F; A))$ which is isomorphic to the invariant subgroup $H^q (F; A) ^G$ where $G=\pi _1 (B, b_0)$. For each $r\geq 2$, there is a injective map $E_r ^{0, q} \rightarrowtail E_2 ^{0,q}$ and a projection $H^q (E; A) \twoheadrightarrow E_{\infty} ^{0, q}$. The composition 
\[
H^q (E; A) \twoheadrightarrow E_{\infty} ^{0, q} \rightarrowtail E_2 ^{0, q} \cong H^q (F; A) ^G
\]
is called the vertical edge homomorphism.
Similarly, there is an isomorphism for the horizontal edge $E_2 ^{p, 0} \cong H^p (B; H^0 (F; A)) \cong H^q (B; A)$, and the horizontal edge homomorphism is defined as the composition
\[
H^p (B; A) \cong E_2 ^{p, 0} \twoheadrightarrow E_{\infty} ^{p, 0} \rightarrowtail H^p (E; A).
\]
The horizontal edge homomorphism coincides with the homomorphisms induced by $\pi : E \to B$.

For each $b\in B$, consider the diagram of maps
\[
\xymatrixcolsep{3.5pc}
\xymatrixrowsep{2.5pc}
\xymatrix{  & H^{n+1} (B, b; A)\ar[d]_<<<{\pi^*}\ar[r]^{j^*} & H^{n+1} (B; A)  \\ H^n (F; A) \ar@{-->}[urr]^<<<<<<<<<<<<<<<<{\tau} \ar[r]_-{\delta_A } & H ^{n+1} (E, F; A)  &}
\]
where $\delta_A $ is the connecting homomorphism for the sequence $F \to E\to (E, F)$ with coefficients in $A$, and $j^*$ is the map induced by the inclusion $B \to (B, b)$. 

\begin{definition}[{\cite[page 186]{McCleary-2001}}]\label{def:Transgression}
There is a map 
\[
\tau : \delta _A ^{-1} (\im \ \pi ^* )  \to H^{n+1} (B; A )/ j^* (\ker \pi ^*)
\]
defined by $\tau (x)= j^* (y)$, where $y \in H^{n+1 } (B, b; A)$ such that $\pi^* (y)=\delta _A (x)$. This map is called a \emph{transgression}. A cohomology class $x \in H^n (F; A)$ is  \emph{transgressive} if there is  $y\in H^n (B, b ; A)$ such that $\pi ^* (y) =\delta_A  (x)$. 
\end{definition}

The differential $d_{n+1} : E_{n+1} ^{0, n} \to E_{n+1} ^{n+1, 0}$ between vertical and horizontal edges on the $E_{n+1}$-page coincides with the transgression $\tau$ (see \cite[Theorem 6.8]{McCleary-2001}). We prove the following proposition on the interaction between the transgression and the connecting homomorphism associated to a short exact sequence.

\begin{proposition}\label{pro:Anticommute} Let $0 \to A\to A' \to A'' \to 0$ be a short exact sequence of abelian groups, and let $\delta _F : H^n (F; A'') \to H^{n+1} (F; A)$ and $\delta _B : H^{n+1} (B; A'') \to H^{n+2} (B; A)$ be the connecting homomorphisms for the corresponding long exact cohomology sequences for the spaces $F$ and $B$. If $x\in H^n (F; A'')$ is transgressive, then $\delta _F (x)$ is also transgressive and $\tau _A ( \delta _F (x)) = - \delta _B (\tau _{A''} (x))$. 
\end{proposition}

\begin{proof}
Since the connecting homomorphisms are natural with respect to chain maps, from the definition of the transgression $\tau$, we see that the equality in the proposition will follow once we show that \begin{equation}\label{eqn:Anticommute} 
\delta_A \circ \delta _F=-\delta _{(E, F)} \circ \delta _{A''}.  
\end{equation}
Note that there is a commuting diagram of chain complexes
\[
\xymatrix{ & 0 \ar[d] & 0 \ar[d] & 0 \ar[d] &  \\
0 \ar[r] & C^*  (E, F; A ) \ar[r] \ar[d] & C^* (E; A) \ar[r] \ar[d] & C^* (F; A) \ar[r] \ar[d] & 0 \\
0 \ar[r] & C^*  (E, F; A' ) \ar[r] \ar[d] & C^* (E; A') \ar[r] \ar[d] & C^* (F; A') \ar[r] \ar[d] & 0 \\
0 \ar[r] & C^*  (E, F; A'' ) \ar[r] \ar[d] & C^* (E; A'') \ar[r] \ar[d] & C^* (F; A'') \ar[r] \ar[d] & 0    \\
& 0 & 0 & 0 &}
\]
and $\delta _A$ and $\delta _{A''}$ are the connecting homomorphisms for the top and bottom horizontal sequences, and $\delta _{(E,F)}$ and $\delta _F$ are the connecting homomorphisms for the left and right vertical sequences. The equality in $(\ref{eqn:Anticommute})$ now follows 
from an easy diagram chase. This equality is proved as Proposition 9.20 in \cite{Osborne-2000}.  
\end{proof}

Consider the exact sequence of abelian groups $0 \to \bbZ \xrightarrow{\times 2 } \bbZ \xrightarrow{m_2} \bbF_2 \to 0$, where $m_2$ is the mod-$2$ reduction map. The connecting homomorphism $ H^n (X; \bbF_ 2 ) \to H^{n+1} (X; \bbZ)$ induced by this short exact sequence, is called the \emph{integral Bockstein operator}, denoted by $\beta_0$. Note that the composition of $\beta _0$ with the induced map $m_2 : H^{n+1} (X; \bbZ ) \to H^{n+1} (X; \bbF_2)$ gives the Bockstein operator $\beta : H^n (X; \bbF_2 ) \to H^{n+1} ( X; \bbF_2)$. The following is a direct consequence of Proposition \ref{pro:Anticommute}.

\begin{corollary}\label{cor:Commute}
Let $F \to E \to B$ be a fibration with path-connected base and fiber spaces, and let $\tau _{\bbZ}$ and $\tau _{\bbF_2}$ denote the transgression maps in the Serre spectral sequences for the fibration with coefficients in $\bbZ$ and $\bbF_2$. Then if $x\in H^n (F; \bbF_ 2)$ is transgressive,  $\beta_0 (x) \in H^{n+1} (F; \bbZ)$ is also transgressive and $\tau_{\bbZ} ( \beta _0 (x)) = \beta_0 (\tau _{\bbF_ 2} (x))$. 
\end{corollary}

\begin{remark}  
Corollary \ref{cor:Commute} is well known to experts in the field. The proof given above provides a direct derivation using standard techniques from homological algebra.
\end{remark}


\section{Integral cohomology of products of real projective spaces}\label{sect:IntegralCohomology} 

Let $G=(\bbZ/2)^r$ and $X$ be a finite-dimensional free $G$-CW-complex homotopy equivalent to $\prod _{i=1} ^k \mathbb{R}P ^{n_i}$.
In this section we describe the integral cohomology rings of $BG$ and $X$.
To avoid the special cases, throughout the section we assume that $n_i\geq 2$ for all $i$. The general case will be handled in the proof of Theorem \ref{thm:Main}.

By K\" unneth's theorem,  there is an isomorphism
$H^* (G; \bbF_2)  \cong \bbF_2 [x_1, \dots, x_r]$ where
$|x_i|=1$ for all $i=1, \dots, r$. For the integral cohomology ring of $G$, first observe that, by K\" unneth's theorem again, the integral cohomology groups of $G$ have exponent 2 in positive dimensions. Using the long exact cohomology sequence for the sequence $0\to \bbZ \xrightarrow{\times 2} \bbZ \xrightarrow{q_2} \bbZ/ 2 \to 0$ of coefficients, we see that for each $n\geq 1$, there is a short exact sequence
\[
0 \to H^n (G; \bbZ) \xrightarrow{m_2} H^n (G;\bbF_2 ) \xrightarrow{\beta_0} H^{n+1} (G; \bbZ) \to 0,
\]
where $m_2: H^* (G; \bbZ) \to H^* (G; \bbF_2)$ is the map induced by the mod-$2$ reduction map $q_2: \bbZ\to \bbZ/2$. Note that since $m_2 : H^{n+1} (G; \bbZ) \to H^{n+1} (G; \bbF_2)$ is injective, the condition $\beta_0(x)=0$ is equivalent to the condition that $\beta(x)=0$.
Hence, via the injective mod-$2$ reduction map,
the positive-dimensional part of $H^*(BG;\mathbb Z)$
identifies with the subring
of $H^* (G; \bbF_2)$ consisting of elements $u \in H^* (G; \bbF_2)$ such that $\beta(u)=0$.
 
To describe the integral cohomology of $G$ as a ring, consider, for each nonempty subset $I \subseteq \{ 1, \dots, r\}$, 
the cohomology class $u_I \in H^* (G; \bbZ)$
defined by 
\[
u_I = \beta_0 \Bigl ( \prod _{i\in I} x_i \Bigr ). 
\]
For $I, J \subset \{ 1,\dots, r\}$, let $I \Delta J$ denote the symmetric difference $(I\backslash J) \cup (J \backslash I )$. 

\begin{proposition}[{\cite[Proposition 6.1]{Benson-Carlson-1992}}]\label{pro:IntCoh}
Let $G\cong (\bbZ/2)^r$. Then $H^* (G; \bbZ)$ is generated by the elements $u_I$, where $\emptyset \neq I \subseteq \{ 1, \dots, r\}$, subject to the relations $2u_{I}=0$ and
\[
u_I u_J =\Bigl (\prod _{j \in I \cap J } u_j \Bigr ) \Bigl ( \sum _{j\in I \cap J } u_{(I \Delta J) \cup \{ j\} } + \sum _{j \in J \backslash I } u_j u _{(I \Delta J) \backslash \{ j\} } \Bigr )
\]
with no restriction on $I$ and $J$ except the assumption that $I$ and $J$ are non-empty. 
\end{proposition}

The integral cohomology of $X\simeq \prod_{i=1} ^k \bbR P ^{n_i}$ can also be calculated using the long exact cohomology sequence for the short exact sequence $0 \to \bbZ \xrightarrow{\times 2} \bbZ \xrightarrow{m_2} \bbZ/2 \to 0$:
\begin{align*}
\cdots \to H^{n-1} (X; \bbF_2) & \xrightarrow{\beta_0} H^n (X; \bbZ) \xrightarrow{\times 2} H^n (X; \bbZ) \xrightarrow{m_2 } H^n (X; \bbF_2 )  \\
& \xrightarrow{\beta_0 } H^{n+1} (X; \bbZ) \to \cdots
\end{align*}

For the cohomology ring of a real projective space, we have
\[
H^* (\bbR P ^n ; \bbZ ) \cong \begin{cases} \bbZ [ s, v] / (2s, s^{m+1}, sv, v^2) &  \text{ if } n=2m+1 \\ \bbZ [s] /  (2s, s^{m+1})  & \text{ if } n=2m
\end{cases}
\]
where $|s|=2$ and $|v|=2m+1$ (see \cite[page 222]{Hatcher-2002}). A direct calculation shows that $m_2(s)=t^2$ and $m_2(v)= t ^{2m+1}$, where $H^* (\bbR P^n; \bbF_2) \cong \bbF_2 [t]/(t^{n+1} )$.
For each $n\geq 0$, we have $H^n (X ; \bbZ )\cong F_n \oplus T_n$, where $F_n$ is a free abelian group and $T_n$ is the torsion subgroup. 
By K\" unneth's theorem, the torsion subgroup $T_n$ has exponent $2$. Using the long exact sequence above, we observe that there is a short exact sequence
\begin{equation}\label{eqn:Short}
0 \to (F_n/ 2F_n) \oplus T_n \xrightarrow{m_2} H^n (X; \bbF_2) \xrightarrow{\beta_0 } T_{n+1} \to 0
\end{equation}
of $\bbF_2$-vector spaces.

Consider the classifying map $q: X\to B\pi$, where $\pi=\pi_1 (X, x_0)$. Since we assumed that $n_i\geq 2$ for all $i$, we have $\pi \cong (\bbZ/2)^k$, hence $H^* (B\pi ; \bbF_2)\cong \bbF_2 [t_1, \dots, t_k]$, where $|t_i|=1$. For each $\emptyset \neq I \subseteq \{1, \dots, k\}$, let 
\[
s_I=\beta_0 \Bigl (\prod _{i\in I}  t_i \Bigr ).
\]
The integral cohomology $H^* (B\pi ; \bbZ)$ is generated by $s_I$, where $\emptyset \neq I \subseteq \{ 1, \dots, k\}$, 
subject to the relations $2s_{I}=0$ and
\begin{equation}\label{eqn:Relations}
s_I s_J =\Bigl (\prod _{j \in I \cap J } s_j \Bigr ) \Bigl ( \sum _{j\in I \cap J } s_{(I \Delta J) \cup \{ j\} } + \sum _{j \in J \backslash I } s_j s_{(I \Delta J) \backslash \{ j\} } \Bigr )
\end{equation}
analogous to the relations given for $u_I$ in Proposition \ref{pro:IntCoh}. We have the following observation.

\begin{lemma}\label{lem:Surjective} The induced map $q^*: H^* (B\pi;\bbZ) \to H^* (X; \bbZ)$ is surjective on the torsion part of $H^* (X; \bbZ)$.
\end{lemma}

\begin{proof} For $n\geq 1$, consider the commuting diagram
\[
\xymatrix{0 \ar[r] & (F_n/2F_n)\oplus T_n \ar[r]^{m_2}  & H^n (X; \bbF_2) \ar[r]^{\beta_0}   & T_{n+1} \ar[r] & 0 \\
0 \ar[r] & H^n (B\pi ; \bbZ) \ar[r]^{m_2} \ar[u]_{q^*} & H^n ( B \pi ; \bbF_2) \ar[r]^{\beta_0} \ar[u]_{q^*}   & H^{n+1} (B \pi; \bbZ) \ar[r] \ar[u]_{q^*} & 0.}
\]

By K\" unneth's theorem $H^* (X; \bbF_2) \cong \bbF_2[t_1, \dots, t_k]/(t_1^{n_1+1}, \dots, t_k ^{n_k+1} )$, hence the middle vertical map $q^*$ is surjective. This implies that the map $q^*$ on the right is also surjective. 
\end{proof}

For each $i\in \{ 1, \dots , k\}$, there is a map $\pi _i : X\to \bbR P^{n_i} $ obtained from the projection map to the $i$-th coordinate. For each $i$ where $n_i$ is odd, let $v_i = \pi _i ^* (v)$ where $v$ is the generator of $H^{n_i} (\bbR P^{n_i}; \bbZ)$. For every $i$ where $n_i$ is even, we take $v_i$ to be zero. By K\"unneth's theorem, the free part of $H^*(X;\mathbb Z)$ is the exterior algebra generated by the classes $\{v_i\}_{i=1} ^k$.
   
\begin{proposition}\label{pro:IntCohX}
    Let $X\simeq \prod_{i=1}^k \bbR P^{n_i}$, $n_i \geq 2$, where $n_i=2m_i$ if $n_i$ is even and $n_i=2m_i+1$ if $n_i$ is odd. 
    Then the integral cohomology of 
    $X$ is generated by the elements $s_I$, where $\emptyset \neq I \subseteq \{ 1, \dots, k\}$, and the elements $v_i$ for all $i$, with degrees $|s_{I}|=|I|+1$ and $|v_i|=n_i$ when $n_i$ is odd, subject to the following relations:
    \begin{enumerate} 
    \item $2s_I =0$ for all $I$,
\item the relations for $s_I$ given in Equation (\ref{eqn:Relations}),
    \item $v_i ^2=0$ for all $i$, and $v_i=0$ when $n_i$ is even,
    \item $v_i s_I=0$ if $i \in I$, and $v_i s_I =s_i ^{m_i} s_{I \cup \{ i\}}$ if $i \not \in  I$, and
    \item $(\prod _{i\in I } s_i ^{m_i}) s_I =0$ for all $I$.
    \end{enumerate}
\end{proposition}

\begin{proof}
The free part of $H^* (X;\bbZ)$  is generated by the set of all $v_i$ where $n_i$ is odd, so the free part is generated by the given generating set.
By Lemma \ref{lem:Surjective}, the induced map $q^* : H^* (B\pi; \bbZ) \to H^* (X; \bbZ)$ is surjective on the torsion part and $H^* (B\pi; \bbZ)$ is generated by the elements $s_I$ where $\emptyset \neq I \subseteq \{1, \dots, k\}$. Hence $H^* (X; \bbZ)$ is generated
by the given set of generators.

The relations given in $(1)$ and $(2)$ come from the relations for the integral cohomology of $B\pi$. Since $v_i =\pi _i ^* (v)$, the equality $v_i ^2=0$ follows from $v^2=0$.
The products appearing in relations $(4)$ and $(5)$ are torsion classes.
They hold because their mod-$2$ reductions hold in \(H^*(X;\mathbb F_2)\).

To complete the proof, we need to show that there are no more relations. Note that the free part of $H^*(X;\mathbb Z)$ is the exterior algebra generated by the classes $v_i$ with $n_i$ odd.
Hence monomials in the $v_i$'s are linearly independent and satisfy no relations other than $v_i^2=0$ and graded commutativity. Since the torsion subgroup intersects the free part trivially, it remains only to analyze relations involving torsion classes. 

Using the relations in $(4)$, any monomial in $s_I$'s and $v_i$'s that has at least some nontrivial $s_I$ term, reduces to an element in the image of $q^* : H^* (B\pi; \bbZ) \to H^* (X;\bbZ)$. Modulo the relations in (2), every element in $H^* (B\pi ; \bbZ)$ can be expressed as a linear combination of elements of the form 
\begin{equation}\label{eqn:Monomial}
(s_{j_1} ^{a_1} \dots s_{j_p} ^{a_p})\cdot s_I,
\end{equation}
where $I =\{ i_1, \dots, i_q\}$ is a nonempty set such that $j_1< \dots< j_p$, $i_1< \dots < i_q$, and $j_p\leq i_q$. The elements of this form are linearly independent in $H^* (B\pi; \bbZ)$
(see the proof of \cite[Proposition 6.1]{Benson-Carlson-1992}). 
Such an element maps to zero under $q^*$ if and only if its mod-$2$ reduction is in the ideal $J=(t_1 ^{n_1+1}, \dots, t_k ^{n_k+1})$. The mod-$2$ reduction of an element in the above form is 
\[
(t_{j_1} ^{2a_1} \dots t_{j_p} ^{2a_p})\cdot \Bigl ( \prod _{i\in I} t_i \Bigr ) \cdot \Bigl ( \sum _{k\in I} t_k \Bigr ).
\]
This polynomial lies in the ideal $J$ if and only if each monomial in the sum has a factor $t_j^{b_j}$ with $b_j \geq n_j+1$. Fix a $k\in I$ and consider the monomial
\[
(t_{j_1} ^{2a_1} \dots t_{j_p} ^{2a_p})\cdot \Bigl ( \prod _{i\in I} t_i \Bigr ) \cdot t_k.
\]
Let $t_j ^{b_j}$ be the factor in this monomial  such that $b_j \geq n_j+1$. Since $b_j\geq 3$, we must have $j=j_l$ for some $l$. 

If $j_l\not \in I$, then $2a_l \geq n_{j_l}+1 \geq 2m_{j_l} +1$. This implies $a_l\geq m_{j_l}+1$. By relation $(5)$ in Proposition \ref{pro:IntCohX} applied to 
$I=\{j_l\}$, we have $s_{j_l} ^{m_{j_l}+1}=0$. This gives $s_{j_l } ^{a_l}=0$, hence the monomial in Equation \ref{eqn:Monomial} is zero. 

Assume that $j_l\in I$ and $j_l\neq k$. Then $2a_l +1 \geq n_{j_l} +1$.  If $n_{j_l}=2m_{j_l}$, this gives $a_l \geq m_{j_l}$. 
As a consequence of relation $(4)$ together with \(v_i=0\) for even \(n_i\),
we have \(s_i^{m_i}s_I=0\) whenever \(n_i=2m_i\) and \(i\in I\).
So, $a_l\geq m_{j_l}$ gives that the monomial in Equation \ref{eqn:Monomial} is zero. If $n_{j_l}=2m_{j_l}+1$, then
$2a_l+1\ge 2m_{j_l}+2$ gives $a_l\ge m_{j_l}+1$. This gives
$s_{j_l}^{a_l}=0$, so the monomial in Equation \ref{eqn:Monomial} is zero.

Finally, assume that $j_l \in I$, and $j_l=k$. Then  $2a_l+2\geq n_{j_l}+1\geq 2m_{j_l}+1$, which implies that $a_l\geq m_{j_l}$. We obtained that for each $k\in I$, either the monomial in Equation \ref{eqn:Monomial} is already zero modulo the stated
relations, or there exists an index $l$ with $j_l=k$ such that $a_l\ge m_k$. If the monomial is not already zero, then this gives that the monomial is divisible by
\[
\Bigl ( \prod_{k\in I}s_k^{m_k} \Bigr ) s_I.
\]
Since this factor is zero by the relation $(5)$ in Proposition \ref{pro:IntCohX}, the monomial in Equation \ref{eqn:Monomial} is zero.  
Hence an element of the form $q^*(u)$ is zero in $H^*(X;\mathbb Z)$ if and
only if $u$ lies in the ideal generated by the relations stated in Proposition \ref{pro:IntCohX}. This completes the proof.
\end{proof}


\section{Free actions with trivial action on integral cohomology}\label{sect:TrivialAction} 

Let $G=(\bbZ/2)^r$ and $X$ be a finite-dimensional free $G$-CW-complex homotopy equivalent to $\prod _{i=1} ^k \mathbb{R}P ^{n_i}$ with $n_i\geq 2$ for all $i$.
Assume that the induced $G$-action on integral cohomology of $X$ is trivial. Consider the Serre spectral sequence for the Borel fibration with mod-$2$ coefficients. Let $t_1,\dots, t_k$ be a basis for $H^1 (X; \bbF_2)$ such that $t_i$ corresponds the $i$-th term in the product, and let $\alpha_1, \dots, \alpha_k \in H^2 (G; \bbF_2)$ be the corresponding k-invariants of the action. By Lemma \ref{lem:Properties}, for every $i$ where $n_i$ is even, $\alpha _i=0$. 

Let $V$ be the vector space generated by $\alpha_1,\dots, \alpha_k$, and $s=\dim V$. Since $r=0$ case is trivial, we have $r\geq 1$. By Proposition \ref{pro:NoCommonZeros}, $\alpha_1,\dots, \alpha_k$ have no nontrivial common zeros in $(\bbF_2)^r$, so we must have $s\geq 1$. By changing the indexing, we can assume that $\alpha_1, \dots, \alpha_s$ is a basis for $V$. 

\begin{definition}\label{def:Irreducible}
We say that the basis $\{\alpha_1, \dots, \alpha_s\}$ for $V$ is \emph{irreducible} if $G$ has no subgroup $H$ of codimension $n\geq 1$ such that the $\bbF_2$-vector space generated by $\{\res^G _H \alpha_i\}_{i=1} ^s$ has dimension $\leq s-n$.
\end{definition}

We have the following lemma.

\begin{lemma}\label{lem:NoLinearSyzygies} Let $\{ \alpha_1,\dots, \alpha_s\}$ be an irreducible basis for $V$. Then, for every linear homogeneous polynomials $l_1,\dots, l_s$ such that $\sum _{i=1}^s l_i \alpha_i=0$, we have $l_1=\dots=l_s=0$.
\end{lemma}

\begin{proof} Let $l'_1,\dots, l'_t$ be a basis for the vector space generated by $l_1,\dots, l_s$. For each $i=1,\dots, s$, we have $l_i=\sum _{j=1} ^t \gamma_{ij} l_j'$ for some $\gamma _{ij} \in \bbF_2$. Then the equation $\sum _{i=1} ^s l_i \alpha_i=0$ gives 
\[
\sum _{i=1} ^s \Bigl ( \sum _{j=1} ^t \gamma_{ij} l_j' \Bigr ) \alpha_i=\sum _{j=1} ^t l_j ' \Bigl ( \sum _{i=1} ^s \gamma _{ij} \alpha_i \Bigr )=0.
\]
Let $q_j :=\sum _{i=1} ^s \gamma_{ij} \alpha_i$ for $j=1,\dots,t$. Since $l_1', \dots, l'_t$ are linearly independent, they form a regular sequence in $\bbF_2[x_1, \dots, x_r]$. 
This follows from the fact that when $l_j'$ are linearly independent, after a change of coordinates they become $x_1,\dots,x_t$, hence form a regular sequence.
 
Since $l_1', \dots, l'_t$ form a regular sequence in $\bbF_2[x_1, \dots, x_r]$, the equation
 $\sum _{j=1} ^t l_j' q_j=0$ implies that $q_j \in (l_1',\dots, l_t')$ for every $j$. 
Let $H$ be the subspace of $G$ defined as $\cap _{i=1} ^ t \ker(l_i')$. Then for every $j$, we have
\[
0=\res^G _H q_j=\sum _{i=1} ^s \gamma _{ij} \res^G _H \alpha_i.
\]

Since the matrix $(\gamma_{ij})$ has rank $t$, these are $t$ linearly independent relations among $\res^G _H \alpha_i$. Therefore the vector space generated by $\res^G_H \alpha_1, \dots, \res^G _H \alpha_s$ has dimension $\leq s-t$. Since $H$ has codimension $t$ in $G$, this contradicts the assumption that the basis $\{\alpha_1,\dots, \alpha_s\}$ is irreducible unless $t=0$. So, we must have $t=0$, hence $l_i=0$ for all $i=1,\dots, s$.
\end{proof}

Since $t_1,\dots,t_k$ is the basis for $H^1(X; \bbF_2)$. For each $i=s+1,\dots,k$, there exist $\lambda_{ij}\in\mathbb F_2$ such that $
\alpha_i=\sum_{j=1}^s \lambda_{ij}\alpha_j.$
Define
\[
t_i'=t_i+\sum_{j=1}^s\lambda_{ij}t_j.
\]
Then \(d_2(t_i')=0\) for every $i\geq s+1$.

\begin{definition} Let $\{t_1, \dots, t_s, t_{s+1}',\dots, t_k'\}$ be a basis of $H^1(X; \bbF_2)$ such that $\{\alpha_i=d_2(t_i)\}_{i=1} ^s$ is a 
basis for $V$, and $d_2(t_i')=0$ for $i=s+1,\dots, k$. Then $\{t_1, \dots, t_s, t_{s+1}', \dots, t_k '\}$ is called a \emph{normalized basis}. If $\{\alpha_1,\dots, \alpha_s\}$ is an irreducible basis, then we say that the normalized basis 
$\{t_1,\dots, t_s, t_{s+1}',\dots, t_k'\}$ is irreducible.  
\end{definition}

We are ready to prove our main result in this section.

\begin{lemma}\label{lem:Square} Let $G=(\bbZ/2)^r$ and let $X$ be a finite-dimensional free $G$-CW-complex homotopy equivalent to $\prod _{i=1} ^k \mathbb{R}P ^{n_i}$, $n_i \geq 2$, such that the induced $G$-action on the integral cohomology is trivial.  Assume that there is a normalized basis 
\[
\{ t_1,\dots, t_s, t_{s+1}',\dots, t_k'\}
\]
for $H^1(X; \bbF_2)$ which is irreducible. Then, for every $i\in \{1,\dots, s\}$, $n_i \equiv 3 \pmod 4$.
\end{lemma}

\begin{proof} By Lemma \ref{lem:Properties}, we already know that $n_i$ is odd for every $i\leq s$. So it remains to rule out the case $n_i\equiv 1 \pmod 4$. Let $i\in\{1,\dots, s\}$ be a fixed index such that $n_i=2m_i+1$ where $m_i=2k_i\geq 2$. We claim that in this case $\beta (\alpha_i)=0$. 

Consider the Serre spectral sequence $E_r ^{p,q} (\bbZ)$ for the Borel fibration with integer coefficients. Since $d_2 (t_i) =\alpha_i$ is a transgression, by Corollary \ref{cor:Commute}, $s_i:=\beta_0 (t_i)$ in $E_3^{0,2} (\bbZ)$ is transgressive and 
\[
d_3 (s_i )=\beta_0 (d_2(t_i))= \beta_0 (\alpha_i).
\]
Since $s_i ^{m_i+1}=0$, we have 
\[
0=d_3 (s_i ^{m_i+1} )= (m_i+1) d_3 (s_i ) s_i ^{m_i}=\beta_0 (\alpha _i) s_i ^{m_i} =\beta(\alpha_i)s_i ^{m_i}
\]
in $E_3 ^{3, 2m_i} (\bbZ)$. Note that multiplication of an integral class $u$ in $H^*(G;\bbZ)$ with a torsion class $s_I$ gives $m_2(u)s_I$. This is why $\beta_0(\alpha_i)$ becomes $\beta(\alpha_i)$.  

Since $\beta(\alpha_i)s_i ^{m_i}=0$ in $E_3 ^{3, 2m_i}(\bbZ)$, there is an element $\theta _i \in E_2 ^{1, 2m_i+1} (\bbZ)$ such that
 $d_2(\theta_i)= \beta (\alpha_i)s_i ^{m_i}$. We have
\[
E_2 ^{1, 2m_i+1} (\bbZ)\cong H^1 (G; H^{2m_i+1} (X; \bbZ)_f) \oplus  H^1 (G; H^{2m_i+1} (X; \bbZ)_t ).
\]
Since $H^1 (G; \bbZ)=0$ and $G$ acts trivially on the integral cohomology of $X$, we have $H^1 (G; H^{2m_i+1} (X; \bbZ) _f )=0$. Thus 
\[
E_2 ^{1, 2m_i+1} (\bbZ)\cong  H^1 (G; H^{2m_i+1} (X; \bbZ)_t )\cong H^1(G;\bbF_2) \otimes H^{2m_i+1} (X; \bbZ)_t.
\]

By Proposition \ref{pro:IntCohX}, every torsion class in $H^{2m_i +1} (X;\mathbb Z)$
is a linear combination of monomials of the form 
\[
s_{j_1}^{b_1}\cdots s_{j_p}^{b_p}s_I
\]
such that $2b_1+\dots+2b_p+|I|+1=2m_i+1$. Since the total degree $2m_i+1$ is odd, $|I|$ must be even. Note that $d_2(s_i)=0$ for all $i$, and $d_2(s_I)=\sum _{i\in I} \alpha_i s_{I\backslash \{ i\}}$.

The only monomials of the above form whose image under $d_2$ involves $s_i ^{m_i}$ is $s_i ^{m_i-1} s_{ij}$ for some $j\neq i$. We can write $\theta_i$ as a sum $u+\theta_i'$ where
\[
\theta_i'= \sum _{p, q \neq i} h_{pq} s_p s_i^{m_i-2} s_{iq}+ \sum _{p<q, i\not \in \{p, q\} } l_{pq} s_i ^{m_i-1} s_{pq} + \sum _{j\neq i} l_j s_i ^{m_i-1} s_{ij}
\]
and $u$ involves only the monomials having $s_i^a$ with $a\leq m_i-3$, or monomials involving $s_i ^{m_i-2} s_{pq}$ with $i\not \in \{p, q\}$, or monomials involving $s_I$ with $|I|\geq 4$. Note that under these conditions $d_2(u)$ does not have any summand that involves $s_i ^{m_i-1}$ or $s_i^{m_i}$. 

Let $y_i=t_i$ for $i\leq s$ and $y_i=t'_i$ for $i\geq s+1$. For each $i,j$, define $w_{ij}:=\beta_0 (y_iy_j)$, and for each $i$, let $w_i:=\beta_0(y_i)$. 
Note that for $i\leq s$, $w_i=s_i$ and for every $i, j\leq s$, $w_{ij}=s_{ij}$. 
For each $i$, let $\alpha_i'=d_2(y_i)$. Then for every $i,j$, we have $d_2(w_{ij})=\alpha_i' w_j +\alpha_j' w_i$.

Since $\{y_1, \dots, y_k\}$ is a basis for $H^1 (X; \bbF_2)$, $\{w_1,\dots, w_k \}$ is a basis for $H^2(X; \bbZ)$. 
The change of basis from $\{t_1,\dots,t_k\}$ to $\{y_1,\dots,y_k\}$ induces an invertible change of basis between the degree-two products $\{t_it_j\}_{i<j}$ and $\{y_iy_j\}_{i<j}$. Applying the integral Bockstein $\beta_0$, which sends $t_it_j$ to $s_{ij}$ and $y_iy_j$ to $w_{ij}$, gives an invertible change of basis between $\{s_{ij}\}_{i<j}$ and $\{w_{ij}\}_{i<j}$. Hence $\{w_{ij}\}_{i<j}$ is an $\bbF_2$-basis for $H^3(X;\mathbb Z)_t$.

Using this basis, the equation for $\theta_i'$ can be written as 
\[
\theta_i'= u'+\sum _{p, q \neq i} h'_{pq} w_p s_i^{m_i-2} w_{iq}+ \sum _{p<q, i\not \in \{p, q\} } l'_{pq} s_i ^{m_i-1} w_{pq} + \sum _{j\neq i} l'_j s_i ^{m_i-1} w_{ij}
\]
for some linear polynomials $h'_{pq}$, $l'_{pq}$ and $l_j'$, and $u'$ is a summand where $d_2(u')$ has no summand involving either $s_i^{m_i}$
or $w_js_i^{m_i-1}$ for any $j$. The coefficients in this new expression come from regrouping the $w_{ij}$-terms, so $l_j'$ may involve coefficients $l_{pq}$ from the original sum. Comparing the
coefficients of \(s_i^{m_i}\) and \(s_i^{m_i-1}\) in $d_2(\theta_i')$, we obtain 
$\beta(\alpha_i)=\sum _{j\neq i} l'_{j} \alpha'_j$, and 
\[ \sum _{p, q\neq i} h'_{p,q} \alpha'_q w_p  s_i ^{m_i-1}+\sum _{p<q,  i\not \in \{p, q\}} l'_{pq} (\alpha'_p w_q +\alpha'_q w_p) s_i ^{m_i-1}  + \sum _{j\neq i} l'_j \alpha'_i w_j   s_i ^{m_i-1}=0.
\]
Note that in this calculation we use the fact that $w_i=s_i$. 

The equation above is of the form $\sum_{j\neq i} A_j w_js_i ^{m_i-1}=0$ for some $A_j\in H^3(G;\mathbb F_2)$.  
Since $m_i=2k_i\ge 2$, Proposition \(\ref{pro:IntCohX}\) implies that $s_i^{m_i-1}\neq 0$. Moreover, for each $j$, the class $s_j s_i^{m_i-1}$ is an admissible monomial, and these monomials are distinct for different values of $j$. Since admissible monomials are linearly independent by the proof of Proposition \ref{pro:IntCohX}, it follows that the classes 
$\{s_j s_i^{m_i-1}\}_{j=1}^k$ are linearly independent.
Since $\{w_j\}_{j=1}^k$ is obtained from $\{s_j\}_{j=1}^k$ by an invertible linear change of basis, the classes
$\{w_js_i^{m_i-1}\}_{j=1}^k$ are also linearly independent.
Hence $A_j=0$ for all $j\neq i$.

From this we obtain that for each $j\neq i$, 
 \[ \Bigl ( \sum _{q\neq i} h'_{j,q} \alpha'_q \Bigr ) + \Bigl ( \sum _{p\neq i, p<j}  l'_{pj} \alpha'_p \Bigr ) + \Bigl ( \sum _{q\neq i, j<q} l' _{jq} \alpha'_q \Bigr ) +l'_j \alpha'_i =0.
\]
Note that since $d_2(t_q')=\alpha'_q=0$ for all  $q\geq  s+1$, the equation above 
involves $\alpha'_q$ only for $q\leq s$, and for every $q\leq s$, $\alpha'_q=\alpha_q$. Hence, by Lemma \ref{lem:NoLinearSyzygies}, all
coefficients of the nonzero classes $\alpha_1,\dots,\alpha_s$ in this equation
must vanish. In particular, $l_j'=0$.
Since this holds for all $j\neq i$ with $j\leq s$, and since $\alpha_j'=0$ for $j\geq s+1$, we conclude that $\beta(\alpha_i)=\sum _{j\neq i} l_j' \alpha'_j=0$. 

To complete the proof, observe that if 
\[
\alpha _i = \sum _j a_j x_j ^2 + \sum _{j<k} a_{jk} x_j x_k,
\]
then $\beta (\alpha _i )= \sum _{j<k} a_{jk} \overline u_{jk}=0$ where 
$\overline u_{jk}:=\beta(x_j x_k)$ for $j<k$. This gives  $a_{jk}=0$ for all $j<k$. Hence $\alpha _i= \ell^2$ where $\ell=\sum _j a_j x_j$. 
If $\ell=0$, then $\alpha_i=0$, contradicting irreducibility. Thus
$\ell\neq 0$, and we may take $H=\ker(\ell)$. Then $\res^G _H \alpha_i=0$. This gives that the dimension of the vector space generated by
$\res _H^G\alpha_1,\dots, \res_H^G\alpha_s$ has dimension $\leq s-1$.
This contradicts the assumption that the basis $\{\alpha_1,\dots,\alpha_s\}$ is irreducible. Therefore no such $i$ can exist. Hence $n_i\equiv 3 \pmod 4$ for every $i\in\{1,\dots,s\}$.
\end{proof}

The following example shows that the assumption $n_i \equiv 1 \pmod 4$ is necessary for the conclusion $\beta(\alpha_i)=0$ in the proof of Lemma \ref{lem:Square}.

\begin{example}\label{ex:Trivial}
Let $Q_8$ denote the quaternion group of order 8. Consider the $Q_8$-action on $V=\bbR ^4$ given by the embedding of $Q_8$ in $\mathbb{H} \cong \bbR ^4$ as the group of unit quaternions. 
For every $m\geq 0$, this gives a $Q_8$-action on $Y=S(\oplus _{m+1} V )\cong S^{4m+3}$. Here, and later in the paper, for every $\bbR G$-module $W$, $S(W)$ denotes the unit sphere in $W$ with respect to a $G$-invariant norm. Note that the central element $z\in Q_8$ of order $2$ acts on $Y$ with the antipodal map. We have $X=Y/\langle z \rangle \simeq \bbR P ^{4m+3}$, and $G=Q_8/\langle z \rangle \cong \bbZ/2 \times \bbZ /2$ acts freely on $X\cong \bbR P^{4m+3}$. 

Since we identified $Q_8$ with the set of unit quaternions $\{ \pm 1, \pm i, \pm j, \pm k\}$ in $\mathbb{H}$, the $Q_8$-action on $\bbR^4$ can be described by 
\[
i \to \begin{bmatrix} 0 & -1 & \phantom{-} 0 & \phantom{-}  0 \\ 1 & \phantom{-} 0 & \phantom{-} 0  & \phantom{-} 0 \\  0  & \phantom{-} 0 & \phantom{-} 0 & -1 \\   0& \phantom{-} 0 & \phantom{-} 1 & \phantom{-} 0  \end{bmatrix} \text{ and } j \to \begin{bmatrix} 0 & \phantom{-} 0 & -1 & \phantom{-} 0 \\ 0 & \phantom{-} 0 & \phantom{-} 0 & \phantom{-} 1 \\ 1 & \phantom{-} 0 & \phantom{-} 0 & \phantom{-} 0 \\ 0 & -1 & \phantom{-} 0 & \phantom{-} 0\end{bmatrix}.
\]
Since both of these matrices have determinant one, the $Q_8$-action on $H^{4m+3} (Y; \bbZ)$ is trivial, hence the induced $G$-action on 
$H^{4m+3} (X; \bbZ )$ is trivial.
Let $\alpha=d_2 (t)$ be the k-invariant of this action. Since the restriction of the $G=Z/2\times \bbZ/2$-action of $X=\bbR P^{4m+3}$ to any cyclic subgroup of $G$ gives a free action of the $\bbZ/2$-action on $X$, we conclude that $\alpha=x_1 ^2 +x_1 x_2 +x_2 ^2$, where $x_1, x_2$ are the generators of $H^* (G; \bbF_2)$. Observe that $\alpha$ is not a square of a 1-dimensional class.
\end{example}

Free $\bbZ/2$-actions on $\bbR P^{4m+1}$ with trivial induced actions on integral cohomology are easy to construct as we show in the following example.

\begin{example}\label{ex:TrivialZ/4}
Consider the $\bbZ/4$-action on $W \cong\bbR^2$ given by rotation by $\pi/2$. Then for every $m\geq 1$, $\bbZ/4$ acts freely on $Y= S( \oplus _{m+1} W)\cong S^{2m+1}$. This gives a $\bbZ/2$-action on $X=Y/(\bbZ/2) \cong \bbR P ^{2m+1}$. Since the $\bbZ /4$-action on $W$ is orientation preserving, the $\bbZ/2$-action on the integral cohomology $H^{2m+1} (X, \bbZ)$ is trivial. The k-invariant of this action is $\alpha =x^2$, where $H^* (\bbZ/2 ; \bbF_2)\cong \bbF_2 [x]$.
\end{example}

By taking products of the group actions given above, we can construct a $G=(\bbZ / 2) ^{k+2l}$ action on $X\cong \prod _{i=1} ^k \bbR P^{2m_i+1} \times \prod _{i=1} ^l \bbR P^{4m_i+3}$ with trivial induced action on integral cohomology. So the upper bound given in Theorem \ref{thm:Main} is realizable by some free action.


\section{Free actions with nontrivial action on integral cohomology}\label{sect:Nontrivial}

We start with an example of a free $(\bbZ/2)^r$-action on a product of real projective spaces with nontrivial action on integral cohomology such that the induced action on the mod-$2$ cohomology is trivial.

\begin{example}\label{ex:Nontrivial}
Consider the $D_8$-action on $S^1$ induced by the symmetries of a square. We can explain this action as the $D_8$-action on the unit sphere $S(W)$, where $W$ is the $2$-dimensional irreducible real representation of $D_8$. If we write $D_8=\langle a, b \, |\, a^2=b^2=(ab)^4=1\rangle$, then
\[
\rho_W(a)= \begin{bmatrix} 1 & \phantom{-} 0 \\ 0 & -1 \end{bmatrix}, \quad \rho_W(b)= \begin{bmatrix} 0 & 1 \\ 1 & 0 \end{bmatrix}.
\]
Let $m \geq 0$ and $V=\oplus _{m+1} W$. Consider the $D_8$-action on $Y=S(V) \cong S^{2m+1}$.
The central element $z=(ab)^2$ in $D_8$ acts freely on $Y$ with the antipodal action. 

Let $X=Y/\langle z \rangle \cong \bbR P^{2m+1}$. Consider the induced $G=D_8/\langle z\rangle \cong \bbZ/2 \times \bbZ/2$-action on $X$. 
The $G$-action has no global fixed point, but the elements $\bar a$ and
$\bar b$ each have fixed points. Thus we can conclude that the k-invariant of this action is equal to $\alpha =x_1x_2$, where $H^* (G; \bbF_2)\cong \bbF_2[x_1, x_2]$. Note that the determinants of $\rho _V (a)$ and $\rho_V (b)$ are equal to $(-1)^{m+1}$. Hence, when $m$ is even, the induced $\overline a$ or $\overline b$ action on $H^{2m+1} (X; \bbZ)\cong \bbZ$ is nontrivial. 

Let $H=\langle \overline{ab} \rangle$. There is an isomorphism $H^{2m+1} (X; \bbZ) \cong \widetilde \bbZ _H$ where 
\[
0 \to \bbZ \to \bbZ[G/H] \to \widetilde \bbZ _H \to 0
\]
is exact. The long exact cohomology sequence associated to this sequence gives 
\[
0 \to H^1 (G , \widetilde \bbZ _H ) \xrightarrow{\delta} H^2 (G; \bbZ ) \to H^2 (G ; \bbZ[G/H]) \to \cdots 
\]
where $H^2 (G; \bbZ[G/H])\cong H^2 (H, \bbZ)$ by Shapiro's lemma.
This gives that 
\[
H^1 (G; \widetilde \bbZ _H)\cong \bbF_2,
\]
generated by $u_H:=\delta ^{-1} \beta_0 (x_1+x_2)$.

We have $\beta (x_1x_2)= x_1 x_2 (x_1+x_2)= m_2 (u_{12})$, where $u_{12}\in H^3 (G; \bbZ)$ is the generator introduced in Proposition \ref{pro:IntCoh}. Hence $\beta_0 (\alpha)= u_{12} \neq 0$. 
Assume that $m=2m'$ for some $m' \geq 1$. Then $X\cong \bbR P^{4m'+1}$. By the argument given in the proof of Lemma \ref{lem:Square}, in the integral Serre spectral sequence for the Borel construction, there is $\theta \in E_2 ^{1, 4m'+1}$ such that $d_2 (\theta)= \beta_0 (\alpha) s^{2m'} =u_{12} s^{2m'}$, where $s=\beta_0 (t)$. Since the product $u_{12} s^{2m'}$ is nonzero, the element $\theta$ is nonzero. Note that 
$E_2 ^{1, 4m'+1}\cong H^1 (G ; H^{4m'+1} (X; \bbZ)) \cong H^1 (G; \widetilde \bbZ_H ) \cong \bbF_2$, hence we can conclude that $\theta$ is equal to $u_H$ under the isomorphism above.  

This example shows that there is a $(\bbZ/2)^r$-action on $\bbR P^{4m+1}$ with nontrivial action on the integral cohomology such that the k-invariant $\alpha$ is not a square and $\beta_0 (\alpha)\neq 0$. This example also shows that the differential $d_2 : E_2 ^{1, 4m+1} \to E_2 ^{3, 4m}$ may be nonzero when the action on the integral cohomology is nontrivial.
\end{example}

\begin{example}\label{ex:NontrivialFree}
In the example above the $G=\bbZ/2 \times \bbZ/2$-action on $Y_1=\bbR P ^{4m+1}$ is not free. However by taking the product with another action of $G$ on $Y_2=\bbR P ^{2l+1}$, we can obtain a free action of $G$ on $X=Y_1\times Y_2$.
To construct the $G$-CW-complex $Y_2$, consider the $\bbZ/4$-action on $S^{2l+1}$ defined by $(z_0, \dots, z_l)\to (iz_0, \dots, iz_l)$. This gives a free $\bbZ/2$-action on $\bbR P^{2l+1}$. 

Let $q: G \to \bbZ/2$ denote the quotient map such that $q(\overline a)=q(\overline b)=1$.
Via $q$, there is a induced $G$-action on $Y_2\cong \bbR P ^{2l+1}$ with kernel $\langle \overline{a}\overline{b} \rangle \cong \bbZ/2$. The $G$-action on
$X=Y_1\times Y_2$ is free since $\langle \overline a \overline b \rangle$ acts freely on $Y_1$. This gives a free $G=(\bbZ/2)^2$ action on $X\cong \bbR P^{4m+1}\times \bbR P^{2l+1}$ where the k-invariants are $\alpha_1=x_1x_2$ and $\alpha_2=(x_1+x_2)^2$. Note that as polynomials in $x_1$ and $x_2$, the k-invariants $\alpha_1$ and $\alpha_2$ have no nontrivial common zeros.
\end{example}

\begin{remark}\label{rem:Mistake}
Note that in the example above, $G=\bbZ/2\times \bbZ/2$ acts freely on $X \simeq \bbR P^{4m+1} \times \bbR P^{2l+1}$ with k-invariants $\alpha_1=x_1x_2$ and $\alpha_2=(x_1+x_2)^2$, and we have $\beta (\alpha_1)=x_1x_2 (x_1+x_2)\neq 0$. This  contradicts a conclusion we had in the proof of \cite[Lemma 8.1]{Yalcin-2000}. This lemma is correct only under the assumption that the $G$-action on the integral cohomology of $X$ is trivial. Because of this, the statement of \cite[Theorem 8.3]{Yalcin-2000} holds only under the same assumption. In the statements of 
\cite[Lemma 8.1]{Yalcin-2000} and \cite[Theorem 8.3]{Yalcin-2000}, 
the assumption that $G$ acts trivially on mod-$2$ cohomology should be replaced by the assumption that $G$ acts trivially on the integral cohomology.
\end{remark}

Before we state the main result of this section, we need to prove some lemmas.

\begin{lemma}\label{lem:Needed} Let $G=(\bbZ/2)^r$ and $\alpha \in H^2 (G; \bbF_2)$ be a cohomology class. If $\beta(\alpha)=\gamma \alpha$ for some $\gamma \in H^1(G; \bbF_2)$, then $\alpha = \eta ( \eta +\gamma)$ for some $\eta \in H^1 (G; \bbF_2)$.
\end{lemma}

\begin{proof} 
If $\gamma=0$, then $\beta(\alpha)=0$. Writing $\alpha$ as a quadratic
polynomial, this implies that the mixed terms vanish, hence
$\alpha=\eta^2$ for some $\eta\in H^1(G;\mathbb F_2)$. Thus
$\alpha=\eta(\eta+\gamma)$.

So we can assume $\gamma \neq 0$. Let $H=\ker \gamma$, then $\beta (\res ^G _H \alpha )=
\res ^G _H (\beta (\alpha))= 0$. This implies that  $\res^G _H \alpha = (u')^2$ for some $u' \in H^1 (H; \bbF_2)$. Let $u \in H^1 (G; \bbF_2)$ be a 1-dimensional class such that $\res^G _H u=u'$. Then $\res^G _H (\alpha +u ^2 )=0$ which gives that $\alpha =u^2 + \gamma \eta$ for some $\eta \in H^1 (G; \bbF_2)$.   
Applying the Bockstein operator to both sides of the equality, we obtain 
$\gamma \alpha = \gamma \eta (\gamma +\eta)$. Since $H^* (G; \bbF_2)\cong \bbF_2 [x_1, \dots, x_r]$ is an integral domain and $\gamma \neq 0$, we conclude $\alpha =\eta ( \gamma +\eta)$.
\end{proof}

We also need the following:

\begin{lemma}\label{lem:Decomposition}
Let $X$ be a finite-dimensional $G$-CW-complex homotopy equivalent to $\prod _{i=1} ^k \bbR P^{n_i}$ such that the induced $G$-action on the mod-$2$ cohomology of $X$ is trivial.  Let $H^n (X; \bbZ)= M_f \oplus M_t$ be a decomposition as an abelian group where $M_f$ is a free subgroup and $M_t$ is the torsion subgroup. Then both $M_f$ and $M_t$ are $G$-invariant, and the decomposition $H^n (X; \bbZ)= M_f\oplus M_t$ is a decomposition of $\bbZ G$-modules. 
\end{lemma}

\begin{proof} Let $M=H^n (X; \bbZ)$ as a $\bbZ G$-module, and let $M=M_f\oplus M_t$ be a decomposition as abelian groups.
Since every element $g\in G$ takes a torsion element to a torsion element, $M_t$ is a $\bbZ G$-submodule of $M$.  

Now we show that $M_f$ is $G$-invariant. Let $m\in M_f$ and $g\in G$.
By (\ref{eqn:Short}), there is a short exact sequence
\[
0 \to M/2M \to H^n (X;\bbF_2) \to T_{n+1} \to 0,
\]
so $M/2M$ is a trivial $\bbF_2 G$-module. Since $M=M_f\oplus M_t$ and $M_t$ has exponent $2$, we have $2M=2M_f$. Since $M/2M$ is a trivial $\bbF_2G$-module, we have $gm-m\in 2M=2M_f$.  Since $m\in M_f$, it follows that $gm\in M_f$. Thus $M_f$ is $G$-invariant. 

Since both $M_f$ and $M_t$ are $\bbZ G$-submodules of $M$, the decomposition $M = M_f\oplus M_t$ is a decomposition of  $\bbZ G$-modules. 
\end{proof}

The main result of this section is the following:

\begin{proposition}\label{pro:Main}
Let $G=(\bbZ/2)^r$ and let $X$ be a finite-dimensional free $G$-CW-complex homotopy equivalent to $\prod _{i=1} ^k \mathbb{R}P ^{n_i}$, $n_i \geq 2$, with trivial action on mod-$2$ cohomology.  Assume that there is a normalized basis 
\[
\{ t_1,\dots, t_s, t_{s+1}',\dots, t_k'\}
\]
for $H^1(X; \bbF_2)$ which is irreducible. Then, for every $i\in \{1,\dots, s\}$, $n_i \equiv 3 \pmod 4$.
\end{proposition}

\begin{proof} By Lemma \ref{lem:Properties}, for all $i$ such that $n_i$ is even, $\alpha_i=0$. So for all $i\leq s$, $n_i$ must be odd. Hence it is enough to prove that there is no $i\leq s$ such that $n_i\equiv 1 \pmod 4$. Let $i\leq s$ be a fixed index such that $n_i=2m_i+1$ for some even $m_i\geq 2$. We claim that in this case $\beta (\alpha_i)=\ell_i \alpha_i$ for some $\ell_i\in H^1 (G; \bbF_2)$. 

Consider the Serre spectral sequence 
\[
E_2 ^{p,q} (\bbZ) =H^p (BG; H^q (X;\bbZ )) \Rightarrow H^{p+q} (X_G; \bbZ)
\]
in integral coefficients for the Borel construction $X \to X_G \to BG$. Note that $H^p (BG; H^q (X; \bbZ))$ is isomorphic to the cohomology of the group $G$ with coefficients in the $\bbZ G$-module $H^q (X; \bbZ)$. By Lemma \ref{lem:Decomposition}, for every $q \geq 0$, we have $H^q(X; \bbZ)\cong H^q (X;\bbZ)_f\oplus H^q(X; \bbZ)_t$ as $\bbZ G$-modules and $G$ acts trivially on the torsion part since it has exponent $2$. Hence
\[
E_2 ^{p,q} (\bbZ) \cong H^p (G; H^q (X; \bbZ)_f) \oplus H^p (G; H^q (X; \bbZ) _t)
\]
where the second summand is isomorphic to $H^p (G;\bbF_2) \otimes H^q (X; \bbZ) _t$.  

The integral spectral sequence $E_2 ^{p,q} (\bbZ)$ has a product structure given by the map
\begin{align*}
H^p(G; H^q (X; \bbZ) \otimes H^{p'} (G; H^{q'} (X; \bbZ) ) & \to H^{p+p'} (G; H^q (X ; \bbZ) \otimes H^{q'} (X; \bbZ) )\\
& \to H^{p+p'} (G; H^{q+q'} (X; \bbZ) )
\end{align*}
where the first map is defined by the cup product in the cohomology of $G$ with nontrivial coefficients, and the second map is induced by the cup product in $H^* (X; \bbZ)$. Note that the products of the cohomology classes in the summand $H^* (G; H^* (X; \bbZ)_t))$ 
stays in this summand, and the Leibniz formula for the differentials holds for such products. 

Let $s_I$, $\emptyset \neq I \subseteq \{1, \dots, k\}$, and $v_j$, $j=1, \dots, k$, be the generators of the integral cohomology ring $H^* (X; \bbZ)$ subject to relations given in Proposition \ref{pro:IntCohX}.
Since $\beta_0 (t_i)=s_i$, by Corollary \ref{cor:Commute}, $s_i$ is transgressive and $d_3(s_i)=\beta_0 (\alpha_i)$. 
By the relations in Proposition \ref{pro:IntCohX}, we have $s_i ^{m_i+1}=0$. Since $s_i\in H^2 (X; \bbZ)_t$, by the remark above, we have 
\[
0=d_3 (s_i ^{m_i+1} )=(m_i+1)d_3 (s_i) s_i ^{m_i}=\beta(\alpha_i)s_i ^{m_i}
\]
in $E_3 ^{3, 2m_i}$. Hence there is an element 
$\theta _i \in E_2 ^{1, 2m_i+1} (\bbZ)$ such that
$d_2 (\theta_i ) = \beta(\alpha_i) s_i ^{m_i}.$

Let $m_2: E_2 ^{*,*} (\bbZ) \to E_2 ^{*,*} (\bbF_2)$ be the morphism of spectral sequences induced by the mod-$2$ reduction map. The mod-$2$ reduction map $m_2 : E_2 ^{p,q} (\bbZ) \to E_2 ^{p,q} (\bbF_2)$ is the composition
\begin{align*}
E_2 ^{p,q } (\bbZ) & \cong H^p (G; H^q (X; \bbZ)_f ) \oplus H^p (G; H^q (X; \bbZ)_t)  \xrightarrow{(m_2, m_2)}  \\
& H^p (G; H^q (X; \bbZ )_f \otimes \bbF_2 ) \oplus H^p (G; H^q (X; \bbZ)_t)  \hookrightarrow H^p (G; H^q (X; \bbF_2)).
\end{align*}
The $\bbF_2$-vector space $H^q(X, \bbF_2)$ is generated by monomials $t_1^{a_1}\cdots t_k^{a_k}$ with $\sum _{i=1} ^k a_i=q$. Since $H^* (X, \bbZ)_f$ is generated by $v_1, \dots v_k$, the cohomology group $H^q (X; \bbZ)_f$ is generated by $v_I:=\prod _{j\in I} v_j$ where $I$ is a set of indices $j$ such that $n_j$ is odd and $\sum _{j\in I} n_j=q$. Hence the image of the mod-2 reduction map is generated by monomials $t_I:=\prod _{j\in I} t_j ^{n_j}$ where $I$ is a set of indices with odd $n_j$ such that $\sum _{j\in I} n_j=q$. 

The image of $H^q(X; \bbZ)_t$ under the mod-2 reduction is
generated by mod-2 reductions of monomials of the form $s_{j_1}^{b_1}\dots s_{j_p} ^{b_p} s_{I}$ where $I=\{i_1,\dots, i_q\}$ with $j_1<\cdots<j_p$, $i_1<\cdots <i_q$ and $j_p<i_q$. Note that none of these give a monomial of the form $t_I$ after mod-2 reduction. So, the monomials involved in the mod-2 reductions of $H^q(X; \bbZ)_f$ and $H^q(X; \bbZ)_t$ are distinct.

Since $m_2$ commutes with the differentials, we have 
\[
d_2 (m_2 (\theta_i ))= m_2 (d_2 (\theta_i ))= m_2 ( \beta(\alpha_i) s_i ^{m_i })=\beta(\alpha_i) t_i ^{2m_i}
\]
in $E_2 ^{3, 2m_i} (\bbF_2)$. Let $M:=H^{2m_i +1} (X; \bbZ)$ and $M=M_f\oplus M_t$ be a decomposition of $M$ as $\bbZ G$-modules. We can write $\theta_i = (\theta _i ^f , \theta_i ^t)$, where $\theta_i ^f \in H^1 (G; M_f)$ and $\theta _i ^t \in H^1 (G; M_t)$. 
We have 
\[
d_2 (m_2 (\theta_i))=d_2 (m_2(\theta_i ^f))+d_2 (m_2 (\theta _i ^t))=\beta(\alpha_i) t_i^{2m_i}.
\]
The image of 
\[
m_2 : H^1 (G; M_f) \to H^1 (G; M_f/ 2M_f)\cong H^1 (G; \bbF_2) \otimes M_f/2M_f 
\]
is generated by elements of the form $\ell\otimes y$ where $\ell\in H^1(G; \bbF_2)$ and $y\in M_f/ 2M_f$. Since $M_f/ 2M_f$ is generated by monomials $t_I=\prod _{j\in I} t_j ^{n_j}$, we can conclude that 
\[
m_2 (\theta _i ^f)= \sum_I \ell_I \otimes t_I 
\]
for some $\ell_I \in H^1 (G; \bbF_2)$. Note that the sum above is over all sets of indices $I$ such that for all $j$, $n_j$ is odd, and $\sum _{j\in I} n_j=2m_i+1$.
Applying $d_2$ to this sum, we obtain
\[
d_2 (m_2(\theta_i ^f) )= \sum _I \sum _{j\in I} \ell_I \alpha_j t_j ^{2m_j} \prod _{h\in I, h\neq j} t_h ^{n_h}
\]
where $n_j=2m_j+1$ for all $j$. Since $m_j\geq 1$ for all $j$, the only term that can contribute to the $t_i^{2m_i}$-component has index $I=\{i\}$. 

Note that the sum above for $d_2 (m_2 (\theta_i ^f))$ is the mod-2 reduction of  
\[
z=\sum _I \sum _{j\in I} \ell_I \alpha_j s_j ^{m_j} \Bigl ( \prod_{h\in I, h\neq j} v_h \Bigr ) =\sum _I \sum _{j\in I} \ell _I \alpha_j s_j ^{m_j-1} s_I \Bigl ( \prod_{h\in I, h\neq j} s_h ^{m_h}   \Bigr )
\]
The last equality follows from the relations given in Proposition \ref{pro:IntCohX}. We can write $z=z'+\ell_i \alpha_i s_i ^{m_i}$ where $\ell_i=\ell_{\{i\}}$.
The equation $d_2 (m_2(\theta_i ^f))+d_2 (m_2 (\theta _i ^t))=\beta(\alpha_i) t_i^{2m_i}$ gives
\[
m_2(z'+\ell_i \alpha_i s_i ^{m_i})+m_2 (d_2(\theta_i ^t))=m_2(\beta(\alpha_i) s_i ^{m_i})
\]
Since the mod-$2$ reductions of torsion basis monomials are linearly
independent and span a direct summand of $H^q(X;\mathbb F_2)$, the induced
map on $H^p(G;-)$ is injective on the torsion summand. Hence we can conclude that
\[
z'+d_2(\theta_i ^t) =(\ell_i \alpha_i+\beta(\alpha_i))s_i ^{m_i}.
\]

As we did in the proof of Lemma \ref{lem:Square}, we can write $\theta^t_i$ as a sum $u+\theta_i'$ where
\[
\theta_i'= \sum _{p, q \neq i} h_{pq} s_p s_i^{m_i-2} s_{iq}+ \sum _{p<q, i\not \in \{p, q\} } l_{pq} s_i ^{m_i-1} s_{pq} + \sum _{j\neq i} l_j s_i ^{m_i-1} s_{ij}
\]
where $d_2(u)$ does not have any summand that involves $s_i ^{m_i-1}$ or $s_i^{m_i}$. We claim that \(z'\) 
has no summand involving either $s_i^{m_i}$ or
$w_js_i^{m_i-1}$. The summands of $z$ are indexed by sets
$I$ such that $\sum_{j\in I}n_j=n_i$. If $i\in I$, then we must have $I=\{i\}$, since all $n_j\geq 2$. 
This is precisely the term $\ell_i\alpha_i s_i^{m_i}$, which has been separated off. Hence every
summand of $z'$ has $i\notin I$, and therefore $z'$ contains no
factor involving $s_i$. In particular, $z'$ does not contribute to the
coefficients of $s_i^{m_i}$ or $w_js_i^{m_i-1}$.

Using the argument given in the proof of Lemma \ref{lem:Square}, the coefficient comparison in the $s_i^{m_i}$-terms gives
\[
\ell_i\alpha_i+\beta(\alpha_i)
=
\sum_{j\neq i} l_j'\alpha_j'.
\]
Applying Lemma \ref{lem:NoLinearSyzygies}, and comparing the coefficient of the $s_i^{m_i-1}$-terms, 
we obtain that $l_j'=0$ for all $j\neq i$ such that $j\leq s$. Since $\alpha_j'=0$ for all $j\geq s+1$, we obtain
$\ell_i \alpha_i +\beta(\alpha_i)=0$.

By Lemma \ref{lem:Needed}, this implies that $\alpha_i=\ell_i' (\ell'_i+\ell_i)$ for some linear polynomial $\ell_i'$. 
Since \(\alpha_i\neq 0\), we have \(\ell_i'\neq 0\).
Taking $H=\ker (\ell'_i)$ we see that 
\[
\dim_{\mathbb F_2}
\left\langle
\operatorname{res}_H^G\alpha_1,\dots,
\operatorname{res}_H^G\alpha_s
\right\rangle
\leq s-1.
\]
This contradicts the assumption that the normalized basis is irreducible. Since $n_i$ is already known to be odd for $i\leq s$, it follows that $n_i\equiv 3 \pmod 4$ for every $i\leq s$.
\end{proof}


\section{Proof of the Main Theorem}\label{sect:MainThm}

We are ready to prove Theorem \ref{thm:Main}. We first observe that it is enough to prove the theorem for the case $n_i\geq 2$ for all $i$.

\begin{proposition}\label{pro:Reduction}
Suppose that Theorem \ref{thm:Main} holds for every   $H=(\bbZ/2)^s$ and every finite-dimensional free
$H$-CW-complex
$X\simeq \prod_{i=1}^k \bbR P^{n_i}$ with $n_i\ge 2$ for all $i$, and  with trivial $H$-action on mod-$2$ cohomology.
Then Theorem \ref{thm:Main} holds for every  $G=(\bbZ/2)^r$ and every finite-dimensional free
$G$-CW-complex $X\simeq \prod_{i=1}^{k'} \bbR P^{n_i'}$
with trivial $G$-action on mod-$2$ cohomology.
\end{proposition}

\begin{proof}
Let $G\cong (\bbZ/2)^r$ and $X$ be a finite-dimensional free $G$-CW-complex homotopy equivalent to $\prod_{i=1}^{k'}\mathbb RP^{n_i'}$, with trivial $G$-action on mod-$2$ cohomology. If for all $i$, $n_i'\geq 2$,
then the conclusion follows from the hypothesis. Otherwise, after reordering
the factors, we can write $X\simeq M\times Y$ where
\[
M=(S^1)^l  \quad \text{ and } \quad
Y=\prod_{i=1}^k \mathbb RP^{n_i},
\]
such that $l\geq 1$ is the number of indices with $n_i'=1$, and $n_i\geq 2$
for all $i$. If $k=0$, we take $Y$ to be a point.

Note that $M=(S^1)^l$ is a compact solvmanifold (see \cite{Jo-Lee-2013}) and $Y$ is a finite-dimensional CW-complex with a finite fundamental group. So we can use the argument given by
Jo and Lee in the proof of \cite[Theorem 2.2]{Jo-Lee-2013}.

We have $\pi _1 (X, x_0) \cong \Gamma \times F$, where $\Gamma \cong \bbZ ^l$ and 
$F \cong (\bbZ /2)^{k}$. Since $G$ acts freely on $X$, there is a group extension
\[
1 \to \pi _1 (X, x_0) \to \pi_1 (X/G, \overline x_0 ) \to G \to 1.
\]
Let $\pi=\pi_1 (X/G, \overline x_0)$. 
Since $F$ is the torsion subgroup of $\pi_1(X,x_0)$, it is characteristic in $\pi_1(X, x_0)$, hence it is normal in $\pi$. 
This gives a commuting diagram of groups
\[
\xymatrix{  & F \ar[d] \ar@{=}[r] & F \ar[d] & & \\ 
1 \ar[r] & \pi_1 (X, x_0) \ar[r] \ar[d] & \pi \ar[r] \ar[d] & G \ar[r] \ar@{=}[d] & 1 \\
1 \ar[r] & \Gamma \ar[r] & \pi/F \ar[r] & G \ar[r] & 1.
}
\]
The $2$-rank of a group $K$, denoted by $\rk_2(K)$, is defined to be the largest $r$ such that  $(\bbZ/2)^r \leq K$. By \cite[Theorem 1.3]{Jo-Lee-2013}, we have 
\begin{equation}\label{eqn:FirstInequality}
\rk_2 G \leq \rk_2 (\pi /F) + l.
\end{equation}

The sequence $1\to \Gamma \to \pi/F \to G\to 1$ is associated to a sequence of covering maps:
\[
\xymatrix{\overline X  \simeq \widetilde M \times Y \simeq Y \ar[d]^{/\Gamma} \\ X \simeq M \times Y \ar[d]^{/G} \\ X/G}
\]
where $\overline X$ is the $\Gamma$-cover of $X$, and $\widetilde M \simeq \bbR^l$ is the universal cover of $M$. 
To see these homotopy equivalences, note that under the homotopy equivalence $X\simeq M\times Y$, the subgroup
$\Gamma\leq \pi_1(X)$ corresponds to the subgroup
$\pi_1(M)\times 1\leq \pi_1(M\times Y)$. Hence the associated
$\Gamma$-cover has homotopy type $\widetilde M\times Y\simeq Y$. 

Since $X$ is a finite-dimensional free $G$-CW-complex, $\overline X$ is a finite-dimensional free $(\pi/F)$-CW-complex (see \cite[Statement (N)]{Whitehead-1949}). Let $H \leq \pi/F$ be a subgroup isomorphic to $(\bbZ /2 )^s$, where $s=\rk_2(\pi/F)$. The restriction of the $(\pi/F)$-action on $\overline X$ to $H$ is also a free action. 

Now we will show that the induced $H$-action on mod-$2$ cohomology of $\overline X$ is trivial.
Since $H\cap \Gamma=1$, the composite $H\leq \pi/F\to G$ is injective. Let $p:\overline X\to X$ be the covering map. Since $p$ is
$H$-equivariant (where $H$ acts on $X$ through the injective homomorphism
$H\hookrightarrow G$), the induced map
\[
p^*:H^*(X;\mathbb F_2)\longrightarrow H^*(\overline X;\mathbb F_2)
\]
is $H$-equivariant. Under the homotopy equivalences
\[
X\simeq M\times Y \quad \text{ and } \quad
\overline X\simeq \widetilde M\times Y\simeq Y,
\]
the map $p^*$ identifies with the map induced by the covering $\widetilde M\times Y\to M\times Y$. Since the projection
\[
H^*(M;\mathbb F_2)\otimes H^*(Y;\mathbb F_2)
\longrightarrow
H^*(Y;\mathbb F_2)
\]
is surjective, $p^*$ is also surjective. The $G$-action on
$H^*(X;\mathbb F_2)$ is trivial, hence the induced $H$-action on
$H^*(\overline X;\mathbb F_2)$ is also trivial.
 
If $k=0$, then $\overline X$ is finite-dimensional and contractible.
A nontrivial finite $2$-group cannot act freely on a finite-dimensional
contractible CW-complex. Hence we must have $H=1$, and
therefore $s=0$. Thus $s\leq \sum_{i=1}^k\mu(n_i)$ holds trivially.

If $k>0$, then the hypothesis applied to the free $H$-CW-complex
\(\overline X\simeq \prod_{i=1}^k\mathbb RP^{n_i}\) gives
\[
s\leq \sum_{i=1}^k\mu(n_i).
\]
Combining this with
\eqref{eqn:FirstInequality}, we obtain
\[
r=\rk_2(G)\leq s+l\leq \left(\sum_{i=1}^k\mu(n_i)\right)+l
=\sum_{i=1}^{k'}\mu(n_i').
\]
This completes the proof of the proposition. 
\end{proof}

Now we are ready to prove Theorem \ref{thm:Main}.

\begin{proof}[Proof of Theorem \ref{thm:Main}]
Let $G=(\mathbb{Z}/2)^r$ and $X$ be a finite-dimensional free $G$-CW-complex 
homotopy equivalent to $\prod _{i=1} ^k  \mathbb{R}P ^{n_i}$ 
such that the $G$-action on the mod-$2$ cohomology of $X$ is trivial.
 
Let $\alpha_1, \dots, \alpha_k$ be the k-invariants in the Serre spectral sequence in $\bbF_2$-coefficients. Consider the $\bbF_2$-vector space $V$ generated by the k-invariants, and assume that the dimension of $V$ is $s$. After reindexing the generators $t_i$, we can assume that $\{t_1, \dots, t_s\}$ is such that $\alpha_1,\dots, \alpha_s$ is a basis for $V$. 
We claim that under this setting the inequality 
\[
r \leq \sum _{i=1} ^s \mu (n_i)
\]
holds. Note that this implies the inequality in the conclusion of the theorem.
We prove the claim by induction on $r$. The case $r=0$ is immediate.

By Lemma \ref{lem:Properties}, for each $i$ where $n_i$ is even, we have $\alpha_i=0$. So $n_i$ is odd for all $i\leq s$.  By Proposition \ref{pro:Reduction}, it is enough to prove the theorem under
the assumption that $n_i\geq 2$ for all $i$. Hence, for $i\leq s$,
we can assume that $n_i=2m_i+1$ with some integer $m_i\geq 1$.

If there is a subgroup $H\leq G$ with codimension $n\geq 1$ such that the vector space $W$ generated by $\{\res^G _H \alpha_i\} _{i=1} ^s$ has dimension $\leq s-n$, then we can reduce the $G$-action on $X$ to the subgroup $H$. Then after reindexing, there is an integer $t\leq s-n$ such that $\{\res^G _H \alpha_i\}_{i=1} ^t$ is a basis for $W$. Applying the induction hypothesis to the free $H$-action on $X$ gives
$r-n \leq \sum _{i=1} ^t \mu (n_i)$. From this we obtain
\[
r\leq n+\sum _{i=1} ^t \mu (n_i) \leq \sum _{i=1} ^s \mu(n_i)
\]
since
$\mu(n_i)\geq 1$ for all odd $n_i$. 
Thus, whenever such a subgroup $H$ exists, the theorem follows from the induction hypothesis. Therefore we may assume that the basis $\{\alpha_1,\dots,\alpha_s\}$ is irreducible (see Definition \ref{def:Irreducible}).

If $s=0$, then all $k$-invariants are zero. By Proposition
\(\ref{pro:NoCommonZeros}\), this is possible only when $r=0$, and the
desired inequality holds. Hence assume that $s\geq 1$. By Proposition \ref{pro:Main}, for every $i\leq s$, $n_i \equiv 3 \pmod 4$. 

By Proposition \ref{pro:NoCommonZeros}, $\alpha_1, \dots, \alpha_s$ have no nontrivial common zeros, hence by \cite[Proposition D]{Cusick-1983}, we have $r \leq 2 s$. Since $n_i\equiv 3 \pmod 4$ for all $i\leq s$, this gives the desired inequality $r \leq \sum _{i=1} ^s \mu(n_i)$. Hence the proof is complete.
\end{proof}
 
It is an open question whether the conclusion of Theorem \ref{thm:Main} holds under the weaker assumption that $X$ has the mod-$2$ cohomology of $\prod _{i=1} ^k  \mathbb{R}P ^{n_i}$, i.e. whether Cusick's conjecture holds in the way it is stated. It is also an interesting question whether the conclusion of Theorem \ref{thm:Main} still holds without the assumption that the induced action on the mod-$2$ cohomology is trivial.

\end{document}